# An efficient algorithm for locating and continuing connecting orbits.

J. W. Demmel, L. Dieci,* and M. J. Friedman†
Computer Science Division, University of California,
Berkeley, CA 94720-1776.
School of Mathematics, Georgia Institute of Technology,
Atlanta, GA, 30332-0160.
Department of Mathematical Sciences, University of Alabama in Huntsville,
Huntsville, AL 35899.

**Abstract**

A successive continuation method for locating connecting orbits in parametrized systems of autonomous ODEs was considered in [9]. In this paper we present an improved algorithm for locating and continuing connecting orbits, which includes a new algorithm for the continuation of invariant subspaces. The latter algorithm is of independent interest, and can be used in different contexts than the present one.

## 1 Introduction.

*Homoclinic* and *heteroclinic orbits,* also called *connecting orbits,* are trajectories connecting equilibrium points of a system of autonomous ordinary differential equations. Computation of connecting orbits is becoming increasingly important, both in dynamical systems research, such as understanding chaotic dynamics, and in a variety of applied problems, including wave propagation in combustion models, chemical reactions, neuronal interactions, solitary waves in fluid, solitons in nonlinear optical fiber, and communication processes in living cells, to name a few. The corresponding numerical problem is that of finding solutions $(u(t), \lambda)$ of the system of autonomous ODEs

$$u'(t) - f(u(t), \lambda) = 0, \quad u(\cdot), f(\cdot, \cdot) \in \mathbb{R}^n, \lambda \in \mathbb{R}^{n_\lambda}, \tag{1}$$

$$\lim_{t \to -\infty} u(t) = u_0, \quad \lim_{t \to +\infty} u(t) = u_1. \tag{2}$$

*Supported in part under NSF DMS-9625813.
†Supported in part by NSF DMS-9404912 and by GM 29123.

Most algorithms for the numerical analysis of connecting orbits reduce (1), (2) to a *boundary value problem* on a finite interval using linear or higher order approximations of stable and unstable manifolds near $u_0$ and $u_1$, respectively. See recent papers by Champneys, Kuznetsov, Sandstede [4], by Doedel, Friedman, Kunin [9] and Moore [13] for the history of the question and the bibliography. Note that in the last work an alternative approach was used based on the arclength parametrization, instead of using time $t$ as a parameter.

The algorithms in [4] use a version of Beyn's continuation algorithm based on projection boundary conditions ([1], [2]). They were implemented in a set of routines, `HomCont`, which are currently part of `AUTO97` [7]. `HomCont` has capabilities for detailed bifurcation analysis of homoclinic orbits and some bifurcations of heteroclinic orbits. It has limited capabilities for locating connecting orbits, namely, a simplified version of the algorithm in [9].

The algorithms in [9] have their primary focus on locating connecting orbits and use a modification of a continuation algorithm based on projection boundary conditions (Friedman, Doedel [11]). They were implemented in an experimental code based on `AUTO94` [8].

In order to have a well posed problem, it is necessary for the boundary conditions to be sufficiently smooth with respect to parameters. Both in [4] and [9], the boundary conditions are defined with respect to bases of stable or unstable eigenspaces of $f_u(u_0, \lambda)$ and $f_u(u_1, \lambda)$. The approach in [4] is to compute an orthonormal basis in the appropriate eigenspace, *at each pseudo arclength continuation step,* using a "black box" routine based on the real Schur factorization and then to adapt this basis to be smooth with respect to parameters, using a technique due to Beyn [2, App. C] which amounts to the solution of a linear system of the dimension of the eigenspace in question. The approach in [9] is to compute initially an orthonormal basis in the appropriate eigenspace via the real Schur factorization and then to continue the real Schur factorization equations (as a part of boundary conditions). At the same time precise convergence of the algorithm in [9] is not clear, and it is somewhat cumbersome to use. In some recent work, [6], Dieci and Eirola provide a general differential equations framework for continuation of the block Schur factorization as well as other matrix factorizations. Reference [6] includes a comprehensive set of references for smooth matrix factorization for parameter dependent matrices.

In this paper we present an improved algorithm for locating and continuing connecting orbits, which includes a new algorithm for the continuation of invariant subspaces (CIS). This CIS algorithm is based on iterative refinement techniques originally due to Stewart [14], and later revisited by Demmel [5]. We provide some new twists to these techniques: (i) we justify these iterative refinement techniques using the differential equations which model continuation of block Schur forms, and (ii) we make use of these differential equations to obtain an accurate approximation of the relevant invariant subspace.

In the end, the new algorithm is more efficient than the algorithms in [4] and [9] and is very robust. In particular, it provides several possible safeguards against fast variation of eigenvalues. It has been implemented in an experimental code based on `AUTO97` which is essentially a modification of the HomCont part of AUTO97 to include the algorithm in [9] for locating and continuing connecting orbits and the CIS algorithm, while preserving the bifurcation analysis part of HomCont.

## 2 An improved algorithm for locating and continuing connecting orbits.

Assume, for simplicity of notation, that the fixed points $u_0$ and $u_1$ are hyperbolic, and the eigenvalues of $f_u(u_0,\lambda)$ and $f_u(u_1,\lambda)$, respectively, satisfy

$$\begin{aligned}\mathrm{Re}\,\mu_{0,n} &\leq ... \leq \mathrm{Re}\,\mu_{0,n_0+1} < 0 < \mu_{0,1} < \mathrm{Re}\,\mu_{0,2} \leq ... \leq \mathrm{Re}\,\mu_{0,n_0},\\ \mathrm{Re}\,\mu_{1,1} &\leq ... \leq \mathrm{Re}\,\mu_{1,n_1} < 0 < \mathrm{Re}\,\mu_{1,n_1+1} \leq ... \leq \mathrm{Re}\,\mu_{1,n}.\end{aligned}$$

In this paper, we will assume that the matrices $f_u(u_{0,1},\lambda)$ are smooth functions of $\lambda$ (for $\lambda$ in an appropriate subset of $\mathbb{R}^{n_\lambda}$.

The method extends to the case $\mu_{0,1}=0$, as in [4]. It also extends to the cases of complex and multiple $\mu_{0,1}$ by a simple modification of Step 0, eq. (11) below, of the algorithm (see [9, Section 4.3] for a computational example). The algorithm requires evaluation of various projections associated with the eigenspaces of $f_u(u_0,\lambda)$ and $f_u(u_1,\lambda)$. Initially these projections are constructed using the real Schur factorizations [12]

$$f_u(u_0,\lambda)=Q_0T_0Q_0^T, \quad f_u(u_1,\lambda)=Q_1T_1Q_1^T.$$

The first factorization has been chosen so that the first $n_0$ columns $q_{0,1},...,q_{0,n_0}$ of $Q_0$ form an orthonormal basis of the right invariant subspace $S_0$ of $f_u(u_0,\lambda)$, corresponding to $\mu_{0,1},...,\mu_{0,n_0}$, and the remaining $n-n_0$ columns $q_{0,n_0+1},...,q_{0,n}$ of $Q_0$ form an orthonormal basis of the orthogonal complement $S_0^\perp$. Similarly, the first $n_1$ columns $q_{1,1},...,q_{1,n}$ of $Q_1$ form an orthonormal basis of the right invariant subspace $S_1$ of $f_u(u_1,\lambda)$, corresponding to $\mu_{1,1},...,\mu_{1,n_1}$, and the remaining $n-n_1$ columns $q_{1,n_1+1},...,q_{1,n}$ of $Q_1$ form an orthonormal basis of the orthogonal complement $S_1^\perp$. In the algorithm below the matrices $Q_0(\lambda)$ and $Q_1(\lambda)$ are assumed to be computed at each continuation step by a "black box" routine, described in Section 3, which ensures their continuity.

The approximate finite interval problem is to find a branch of solutions $(u(t),\lambda,u_0,u_1,T)$, $u\in C^1([0,1],\mathbb{R}^n)$, $\lambda\in\mathbb{R}^{n_\lambda}$, provided $n_\lambda=n-(n_0+n_1)+2$; $n_\lambda\geq 0$, $u_0,u_1\in\mathbb{R}^n$, $T>0$, is the length of the time interval, for some small $\epsilon_0,\epsilon_1>0$, of the time-scaled differential equation

$$u'(t)-Tf(u(t),\lambda)=0, \;\; 0<t<1, \tag{3}$$

subject to left boundary conditions

$$\begin{aligned}(u(0)-u_0)\cdot q_{0,n_0+i}(u_0,\lambda) &= 0, \;\; i=1,...,n-n_0, &(4)\\ |u(0)-u_0| &= \epsilon_0, &(5)\end{aligned}$$

right boundary conditions

$$\begin{aligned}(u(1)-u_1)\cdot q_{1,n_1+i}(u_1,\lambda) &= 0, \;\; i=1,...,=n-n_1. &(6)\\ |u(1)-u_1| &= \epsilon_1 &(7)\end{aligned}$$

stationary state conditions

$$f(u_0, \lambda) = 0, \tag{8}$$
$$f(u_1, \lambda) = 0. \tag{9}$$

**Remark 1** *Initially, we perform time integration to obtain a, typically, crude orbit with initial point $u(0) \in S_0$ but the terminal point $u(1) \notin S_1$, in general. Hence $\tau_i$ defined by*

$$\tau_i = (u(1) - u_1) \cdot q_{1,n_1+i}(u_1, \lambda)/\epsilon_1, \quad i = 1, ..., n_\tau = n - n_1, \tag{10}$$

*are, in general, nonzero, and the initial connecting orbit on the branch of connecting orbits is found via a sequence of homotopies that locate successive zero intercepts of the $\tau_j$ in (10). In each homotopy step we compute a branch, i.e., a one-dimensional manifold, of solutions. For this we must have $n_c - n_v = n - 1$, where $n_c$ is the number of constraints, and $n_v$ is the number of scalar variables. We keep $u(1)$ free, while $u\,(0)$ is allowed to vary on the surface of the sphere in $S_0$ of radius $\epsilon_0$ : (i) according to equation (12) below at steps 0 through $n_0$, when $\lambda$ is fixed and $c_i$ in (12) play the role of the control parameters; and (ii) according to equation (5) at steps $n_0 + 1$ through $n_0 + n_\lambda$, when $\lambda$ varies.*

Let $S_{0,k},\ k = 1, ..., n_0$, be the right invariant subspace of $f(u_0, \lambda_0)$ corresponding to the eigenvalues $\mu_{0,1}, ..., \mu_{0,k}$. Then the first $k$ columns $q_{0,1}, ..., q_{0,k}$ of $Q_0$ form an orthonormal basis of $S_{0,k}$ and the remaining $n-k$. columns $q_{0,k+1}, ..., q_{0,n}$ of $Q_0$ form an orthonormal basis of the orthogonal complement $S_{0,k}^{\perp}$.

1. Initialization.
   Step 0. Initialize the problem parameter vector $\lambda$, and set the algorithm parameters $\epsilon_0$ and $T$ to small, positive values, so that $u(t)$ is approximately constant on $[0, T]$. Set
   $$u(t) = u_0 + \epsilon_0 c_1 q_{01}, \quad 0 \le t \le 1, \tag{11}$$
   or
   $$u(t) = u_0 + \epsilon_0 c_1 q_{01} e^{\operatorname{Re}\mu_{0,1} t}, \quad 0 \le t \le 1, \quad \operatorname{Re}\mu_{0,1} > 0,$$
   $\epsilon_1 = |u(1) - u_1|$, $c_1 = 1$, or $-1$ and $c_2 = ... = c_{n_0} = 0$.

2. Locating a connecting orbit, $\lambda$ is fixed.
   Step 1. Time integration to get an initial orbit. Compute a solution branch to the system (3), (4), (10), (7), and
   $$(u\,(0) - u_0) \cdot q_{0,i}(u_0, \lambda)/\epsilon_0 - c_i = 0, \quad i = 1, ..., n_c = n_0, \tag{12}$$
   in the direction of increasing $T$, until $u(1)$ reaches an $\epsilon_1$-neighborhood of $u_1$, for some $\epsilon_1 > 0$. Scalar variables are $T$, $\epsilon_1 \in \mathbb{R}$, $\tau \in \mathbb{R}^{n-n_1}$. There are $n$ differential equations with $n_c = 2n - n_1 + 1$ constraints and $n_v = n - n_1 + 2$ scalar variables, and hence $n_c - n_v = n - 1$. In practice one typically continues until $\epsilon_1$ stops decreasing, its value being not necessarily small.
   Step $k$, $k = 2, ..., n_0$ (for $n_0 > 1$ ). Compute a branch of solutions to the system (3)-(5), (10), (7), (12) to locate a zero of, say, $\tau_{k-1}$ (while $\tau_1, ..., \tau_{k-2} = 0$, fixed). Free scalar variables are $\epsilon_1, c_1, ..., c_k, \tau_{k-1}, ..., \tau_{n-n_1}$. There are $n$ differential equations with $n_c = 2n - n_1 + 2$ constraints and $n_v = n - n_1 + 3$ scalar variables, and hence $n_c - n_v = n - 1$.

3. Locating a connecting orbit, $\lambda$ varies.
Step $k$, $k = n_0 + 1, ..., n_0 + n_\lambda \equiv n - n_1 + 1$. Compute a branch of solutions to the system (3)-(5), (10), (7)-(9) to locate a zero of, say, $\tau_{k-1}$(while $\tau_1, ..., \tau_{k-2} = 0$, fixed). Free scalar variables are $\epsilon_1, \tau_{k-1}, ..., \tau_{n-n_1}, \lambda_1, ..., \lambda_{k-n_0} \in \mathbb{R}$, $u_0, u_1 \in \mathbb{R}^n$. There are $n$ differential equations with $n_c = 4n - n_0 - n_1 + 2$ constraints and $n_v = 3n - n_0 - n_1 + 3$ scalar variables, and hence $n_c - n_v = n - 1$.

4. Increasing the accuracy of the connecting orbit.
Compute a branch of solutions to the system (3)-(9) in the direction of decreasing $\epsilon_1$ until it is 'small'. Free scalar variables are $\epsilon_1, T, \lambda_1, ..., \lambda_{n_\lambda - 1} \in \mathbb{R}$, $u_0, u_1 \in \mathbb{R}^n$. As before, $n_c = 4n - n_0 - n_1 + 2$, $n_v = 3n - n_0 - n_1 + 3$.

5. Continue the connecting orbit.
Compute a branch of solutions to the system (3)-(9). Free variables are the (real) scalar $T$, and $\lambda_1, \ldots, \lambda_{n_\lambda}$ , and the vectors $u_0, u_1 \in \mathbb{R}^n$. As before, $n_c = 4n - n_0 - n_1 + 2$, $n_v = 3n - n_0 - n_1 + 3$. Alternatively, a phase condition

$$\int_0^1 (u^{'}(t) - q^{'}(t)) \cdot u^{''}(t) \; dt = 0 \qquad (13)$$

may be added if $T$ is kept fixed and $\epsilon_0$ and $\epsilon_1$ are allowed to vary. Here $q(t)$ is a previously computed orbit on the branch.

**Remark 2** *In principle, our algorithm of continuation of invariant subspaces of $f_u(u_0, \lambda)$ and $f_u(u_1, \lambda)$ in Section 3 breaks down only if two eigenvalues, one associated with the subspace being continued and one not, approach the same point on the imaginary axis (one from the left and one from the right). In this case the algorithm should stop anyway since a bifurcation is being approached.*

# 3 Continuation of Invariant Subspaces.

Let $A(\lambda) \in \mathbb{R}^{n \times n}$ denote one of the following: $f_u(u_0, \lambda)$, $f_u(u_1, \lambda)$. The basic continuation algorithm requires at each pseudo arclength continuation step computation of a right invariant (typically, stable or unstable) $m$-dimensional subspace $S(\lambda)$ of $A(\lambda)$. In general, the function $A$ is smooth in $\lambda$ (say, differentiable), and it is important that $S(\lambda)$ be also smooth, as otherwise convergence difficulties can be expected.

In this section we show how to constructively obtain smooth and orthogonal $Q(\lambda) = [Q_1(\lambda)|Q_2(\lambda)] \in \mathbb{R}^{n \times n}$, $Q_1(\lambda) \in \mathbb{R}^{n \times m}$, $Q_2(\lambda) \in \mathbb{R}^{n \times (n-m)}$, so that $Q_1(\lambda)$ span $S(\lambda)$, and $Q_2(\lambda)$ span the orthogonal complement $S(\lambda)^\perp$. Moreover, if we let $S_k(\lambda)$, $k = 1, ..., m$, be the right invariant subspaces of $A(\lambda)$ corresponding to the first $k$ eigenvalues $\mu_1, ..., \mu_k$ of $A$, then the first $k$ columns $q_1, ..., q_k$ of $Q(\lambda)$ form an orthonormal basis of $S_k(\lambda)$, and the remaining $n - k$ columns form an orthonormal basis of the orthogonal complement $S_k(\lambda)^\perp$.

To justify our construction, we should recall that in our continuation procedure we parametrize a solution branch in terms of so called pseudo arclength; let $s$ denote the pseudo arclength variable [1]. Thus, both fixed points $u_0$ and $u_1$, see (8)–(9), as well as the parameter(s) $\lambda$ are smooth functions of $s$. The matrix valued function $A: \lambda \in \mathbb{R}^{n_\lambda} \to \mathbb{R}^{n\times n}$ can thus be viewed as a smooth function from $s \in \mathbb{R} \to \mathbb{R}^{n\times n}$. As a consequence, we can and will think of invariant subspaces' continuation with respect to the scalar pseudo arclength variable $s$. For this reason, we will abuse of notation and freely write $A(s)$ for $A(\lambda)$. In what follows, a "˙ " will indicate differentiation with respect to $s$.

**Remark 3** *If $n_\lambda = 1$ in (1), then it is conceivable that one may be able to perform continuation with respect to $\lambda$ itself, rather than reparametrizing the problem by arclength.*

The basic issue is the following: "Suppose that initially we have the (real) block Schur factorization

$$A(0) = Q(0)T(0)Q^T(0), \quad Q(0) = [Q_1(0)|Q_2(0)], \tag{14}$$

where $T(0)$ is block upper triangular ($T_{ii}(0),\ i = 1, 2$, are not required to be triangular)

$$T(0) = \begin{bmatrix} T_{11}(0) & T_{12}(0) \\ 0 & T_{22}(0) \end{bmatrix},$$

the columns of $Q_1(0)$ span an invariant subspace $S(0)$ of $A(0)$, and the columns of $Q_2(0)$ span the orthogonal complement $S(0)^\perp$. We want to obtain a block Schur factorization for the matrix $A(s)$, close to $A(0)$, exploiting (if possible) the work already done to obtain the block Schur form for $A(0)$". This problem fits within the general framework developed in [6]. Suppose that the matrix $A(s)$ has two groups of eigenvalues, $\Lambda_1(s)$ and $\Lambda_2(s)$, which stay disjoint for all $s$ around 0. Then, in an interval about $s = 0$ there is the smooth factorization

$$A(s) = Q(s)T(s)Q^T(s), \quad Q(s) = [Q_1(s)|Q_2(s)], \tag{15}$$

where $T(s)$ is in block Schur form $T(s) = \begin{bmatrix} T_{11}(s) & T_{12}(s) \\ 0 & T_{22}(s) \end{bmatrix}$. Here, $T_{11}$ has eigenvalues $\Lambda_1$ and $T_{22}$ has eigenvalues $\Lambda_2$. A constructive procedure to obtain the factorization $QTQ^T$ of $A$ can be based upon differential equations models. As in [6], we differentiate the relations $A = QTQ^T$ and $Q^T(s)Q(s) = I$, let $H := Q^T\dot{Q}$, and obtain

$$\begin{aligned} \dot{T} &= Q^T\dot{A}Q + TH - HT\,, & (16)\\ \dot{Q} &= QH\,. & (17) \end{aligned}$$

Now we use triangularity of $T$ and the fact that $H$ must be skew symmetric, $H^T = -H$; we partition $H$ in the same way as $T$, and then can determine $H_{12}$ by solving the Sylvester equation

$$T_{22}H_{12}^T - H_{12}^T T_{11} = (Q^T\dot{A}Q)_{12}\,. \tag{18}$$

The blocks $H_{11}$ and $H_{22}$ are not uniquely determined, and we may set them to 0 (in any case, they must be chosen skew symmetric). Thus, in principle, we can solve (16)-(17) subject to initial conditions (ICs) obtained from the factorization of $A(0)$, in order to obtain a smooth path of block Schur factorizations.

[1] Typically, this is the general procedure of most practical continuation algorithms; e.g., this is the strategy implemented in `AUTO`.

We can accumulate the transformations in such a way that we are always looking for corrections close to the identity. To be more precise, we can rewrite

$$Q(s) = Q(0)U(s)\,, \text{ with } U(0) = I\,, \tag{19}$$

and use this in (17), thereby obtaining a differential equation for $U$, $\dot{U} = UH$ (notice that $H$ is the same as before). We now look for exact solutions to the $U$ differential equation. Partition $U(s) = [U_1(s),\ U_2(s)] = \begin{bmatrix} U_{11}(s) & U_{12}(s) \\ U_{21}(s) & U_{22}(s) \end{bmatrix}$, with same block dimensions as $T$. Since $U(0) = I$, there is an open interval about 0 where we can require that $U_1$ has the structure $U_1 = \begin{bmatrix} I \\ U_{21}U_{11}^{-1} \end{bmatrix} U_{11}$. Next, we define

$$Y(s) \ := \ U_{21}(s)U_{11}^{-1}(s)\,, \tag{20}$$

and use the orthogonality relation $U_1^TU_1 = I$ to obtain $U_1 = \begin{bmatrix} I \\ Y \end{bmatrix} (I + Y^TY)^{-1/2}$. In a similar way, for $U_2(s)$ we use $U_2^TU_2 = I$ and $U_1^TU_2 = 0$, to eventually obtain

$$U(s) = \left[ \begin{pmatrix} I \\ Y(s) \end{pmatrix} (I + Y^TY)^{-1/2}, \begin{pmatrix} -Y^T(s) \\ I \end{pmatrix} (I + YY^T)^{-1/2} \right]. \tag{21}$$

Thus, we need to find the matrix $Y \in \mathbb{R}^{(n-m)\times m}$ in (21). Define $E(s)$ by

$$Q^T(0)A(s)Q(0) = Q^T(0)[A(0) + (A(s) - A(0))]Q(0) =: T(0) + E(s) = \begin{bmatrix} \widehat{T}_{11} & \widehat{T}_{12} \\ E_{21} & \widehat{T}_{22} \end{bmatrix}. \tag{22}$$

Now we substitute $Q(s)$ given by (19), (21) and $A(s)$ obtained from (22) into the invariant subspace relation:

$$Q_2^T(s)A(s)Q_1(s) = 0, \tag{23}$$

to obtain the following algebraic Riccati equation for $Y$ :

$$\widehat{T}_{22}Y - Y\widehat{T}_{11} = -E_{21} + Y\widehat{T}_{12}Y, \tag{24}$$

or

$$F(Y) = 0\,, \quad F(Y) := \widehat{T}_{22}Y - Y\widehat{T}_{11} + E_{21} - Y\widehat{T}_{12}Y. \tag{25}$$

**Remark 4** *If we think of $S(0)$ and $Q(0)$ as approximations to $S(s)$ and $Q(s)$, then we may interpret the form (21) as well as the resulting Riccati equation (24) as an iterative refinement technique to improve the accuracy of computed invariant subspaces. This is the viewpoint in the works of Stewart [14], Dongarra, Moler and Wilkinson [10], Chatelin [3], and Demmel [5]. However, we prefer to think of (21) as the exact solution of the differential equation $\dot{U} = UH$ (see (17) and (19)). We will exploit this fact later.*

**How to solve (24).** The following two iterative methods have been often advocated to solve the Riccati equation (24) (or (25)):

1. The iteration [14]:

$$\widehat{T}_{22}\Delta_k - \Delta_k\widehat{T}_{11} = -F(Y_{k-1}) \ , \ \ Y_k = \Delta_k + Y_{k-1} \, , \tag{26}$$

   with $Y_0 = 0$, $k = 1, 2, \ldots$ .

2. The Newton iteration [5]:

$$(\widehat{T}_{22} - Y_{k-1}\widehat{T}_{12})\Delta_k - \Delta_k(\widehat{T}_{11} + \widehat{T}_{12}Y_{k-1}) = -F(Y_{k-1}) \ , \ \ Y_k = \Delta_k + Y_{k-1} \, , \tag{27}$$

   with $Y_0 = 0$, $k = 1, 2, \ldots$ .

Therefore, we need to solve a Sylvester equation in the inner loop of the iterative refinement. This can be effectively done by using LAPACK routines. As shown in the convergence analysis for the iterations (26) and (27) by Stewart [14] and Demmel [5], respectively, if we let

$$\kappa = \frac{\|\widehat{T}_{12}\|_F \|E_{21}\|_F}{\mathrm{sep}^2(\widehat{T}_{11}, \widehat{T}_{22})}, \tag{28}$$

then under the assumptions $\kappa < 1/4$ and $\kappa < 1/12$, the iterations (26) and (27) converge, linearly and quadratically, respectively. In (28), $\| \cdot \|_F$ is the Frobenius norm of a matrix.

The parameter $\kappa$ can be interpreted as follows. Its numerator, $\|\widehat{T}_{12}\|_F \|E_{21}\|_F$ measures the quality of the initial approximate invariant subspace: it will be small when the approximation is good, and the factor $\|E_{21}\|_F$ will be zero if and only if the initial approximation is in fact correct. The function $\mathrm{sep}(\widehat{T}_{11}, \widehat{T}_{22})$ in the denominator is the smallest singular value of the operator which maps $Y$ to $\widehat{T}_{22}Y - Y\widehat{T}_{11}$, and it measures the separation of the spectra of $\widehat{T}_{11}$ and $\widehat{T}_{22}$. If $\mathrm{sep}(\widehat{T}_{11}, \widehat{T}_{22})$ is small, it means that some eigenvalues of $\widehat{T}_{11}$ and $\widehat{T}_{22}$ can be made to merge with small changes in $\widehat{T}_{ii}$; this means that the invariant subspaces belonging to the two parts of the spectrum are unstable and hard to compute. Thus $\kappa$ will be small if we start with a good initial approximate invariant subspace and if the eigenvalues associated with that subspace are well separated from the remainder of the spectrum. Oversimplifying, both algorithms converge if (i) the spectra of $T_{11}(0)$ and $T_{22}(0)$ are far enough apart, and (ii) the perturbation $E(s)$ of $T(0)$ is small enough.

**Remark 5** *We point to some additional safeguards, which ensure proper performance of the algorithms (26) and (27).*

1. *Note that from (22) we have* $\|E_{21}\|_F \le \|E\|_F \le \|A(s) - A(0)\|_F$. *Hence, from (28), if*

$$\alpha\|A(s) - A(0)\|_F \le \frac{1}{4}, \quad \alpha = \frac{\|\widehat{T}_{12}\|_F}{sep^2(\widehat{T}_{11}, \widehat{T}_{22})}, \tag{29}$$

   *then* $\kappa < 1/4$. *Therefore we are guaranteed that for any matrix* $X \in \mathbb{R}^{n\times n}$ *with* $\alpha\|X - A(0)\|_F \le \frac{1}{4}$ *we can find the invariant subspaces* $S$ *and* $S^\perp$ *by Algorithm (26),* *say. And no eigenvalue of $A\mid_S$ can 'merge' with an eigenvalue from $A\mid_{S^\perp}$, where $A\mid_S$ denotes the restriction of $A$ to $S$, etc.. In other words, the eigenvalues of $A\mid_S$ and $A\mid_{S^\perp}$ remain separated for all matrices in a ball around $A(0)$ of radius $1/(4\alpha)$. Thus, employing the safeguard (29) ensures that the iterations (26) and (27) will diverge only when a small perturbation $A(s)$ of $A(0)$ will make an eigenvalue from $S(s)$ and an eigenvalue from $S(s)^\perp$ coalesce.*

2. *The quantity $\|E_{21}\|_F/sep(\widehat{T}_{11},\widehat{T}_{22})$ can be interpreted as the tangent of the angle between the subspaces spanned by $Q_1(0)$ and $Q_1(s)$. Hence the convergence of our algorithms implies that this angle always stays between 0 and $\pi/2$. If one wants to monitor this angle, a convenient measure is the sine of the angle given [12] by*

$$\mathrm{dist}(Q_1(0),Q_1(s)) = \sqrt{1-\sigma^2_{\min}\left(Q_1^T(0)Q_1(s)\right)}.$$

*Substituting into this equation the formula for $Q_1(s)$ gives*

$$\mathrm{dist}(Q_1(0),Q_1(s)) = \sqrt{1-\sigma^2_{\min}\left(I+Y^TY\right)^{-1/2}} = \frac{\|Y\|_2}{\sqrt{1+\|Y\|_2^2}}\,. \tag{30}$$

The previous discussion has been motivated by the standard linear algebra viewpoint, that is of "how to refine the initial trivial estimate $U=I$". In particular, this gave us the initial conditions $Y_0=0$ for the iterations (26) and (27). However, from the point of view of continuation, this is usually not the best strategy, since it amounts to starting with the old solution as initial guess. We can do better by using the differential equation formulation and taking an Euler approximation to the solution.

Recall that the matrix $Y(s)$, solution of (24), is related to $U$ by (20): $Y(s)=U_{21}(s)U_{11}^{-1}(s)$. Here, $U_1=\begin{bmatrix}U_{11}\\U_{21}\end{bmatrix}$, and $U=[U_1,\,U_2]$ solves the differential equation $\dot U = UH$ (see (17)–(19)). That is

$$\dot U_1 = U_1H_{11}+U_2H_{21},\quad \dot U_2 = -U_1H_{21}^T+U_2H_{22}. \tag{31}$$

Recall that the factor $H_{21}(s)$ is determined by (18), but there is freedom insofar as the choice of $H_{11}$ and $H_{22}$. We now show that $Y(s)$ is unaffected by the choices made for $H_{11}$ (and $H_{22}$). In particular, this implies that the range of $s$–values guaranteeing invertibility of $U_{11}(s)$, and hence the representation (21), is unaffected by the choice of $H_{11}$ and $H_{22}$ in (31), and hence there is little reason not to set $H_{11}$ and $H_{22}$ both equal to 0 when we represent $U$ as in (21).

**Lemma 1** *Let $U(s)=[U_1(s),\ U_2(s)]$ be the solution of (31) with $U(0)=I$, obtained by setting $H_{11}(s)=0$ and $H_{22}(s)=0$ in (31), while determining $H_{21}(s)$ by (18). Let $\tilde U(s)=[\tilde U_1,\ \tilde U_2]$ be the solution of (31) with $\tilde U(0)=I$, and $H_{11}$ and $H_{22}$ nonzero skew symmetric matrices. Then, we have*

$$\tilde U_1(s) \;=\; U_1(s)C(s)\;,$$

*where $C(s)\in\mathbb{R}^{m\times m}$ is orthogonal. Therefore, partitioning $U_1=\begin{bmatrix}U_{11}\\U_{21}\end{bmatrix}$, and similarly for $\tilde U_1$, we have $Y(s)=U_{21}(s)U_{11}^{-1}(s)=\tilde U_{21}(s)\tilde U_{11}^{-1}(s)$.*

**Proof.** Let $p(\lambda, s)$ be the characteristic polynomial of the matrix $A(s)$. Then, we have the factorization $p(\lambda, s) = p_1(\lambda, s)p_2(\lambda, s)$ where the factors $p_1$ and $p_2$ have no common root, and the roots of $p_1$ give $\Lambda_1$, while the roots of $p_2$ give $\Lambda_2$. Thus, we have

$$\begin{aligned} p_1(Q^T(0)A(s)Q(0), s)U_1(s) &= 0 \,, \quad \text{and} \\ p_1(Q^T(0)A(s)Q(0), s)\tilde{U}_1(s) &= 0 \,. \end{aligned}$$

Therefore, since $\text{rank}(p_1(A(s), s)) = n - m$, then $\tilde{U}_1(s) = U_1(s)C(s)$, where $C(s)$ is orthogonal (since $\tilde{U}_1^T(s)\tilde{U}_1(s) = I$). ■

Now, because of Lemma 1 we can just integrate the differential equations (31) with $H_{11}$ and $H_{22}$ equal to 0:

$$[\dot{U}_1, \ \dot{U}_2] = [U_2 H_{21} \,, \ -U_1 H_{21}^T] \,, \quad U(0) = I \,.$$

So, the idea is to approximate the solution of these differential equations from 0 to $s$ in order to get an approximation $V_1$ to $U_1(s)$ and thus obtain an approximation to $Y(s)$ from $V_{21}V_{11}^{-1}$, where we have partitioned $V_1 = \begin{bmatrix} V_{11} \\ V_{21} \end{bmatrix}$. We use a forward Euler step to obtain $V_1$. [2]

1. Initialization.
   Set $U(0) = I$ and obtain $H_{21}(0)$ from solving the Sylvester equation
   $$T_{22}(0)H_{21}(0) - H_{21}(0)T_{11}(0) = -\,(0, \ I)\, Q(0)^T \dot{A}(0)Q(0)\, (I, \ 0)^T \,. \tag{32}$$

2. Euler step.
   Let $V_1 = U_1(0) + sU_2(0)H_{21}(0)$. Obtain initial guess for the Riccati equation (24) (or (25)):
   $$Y_0 = V_{21}V_{11}^{-1}. \tag{33}$$

**Remark 6** *We now look at how this new initial guess $Y_0$ in (33) impacts convergence of the iterations (26) and (27). We also make some comments on expense.*

1. *Observe that the value of $Y_0$ in (33) is nothing but $Y_0 = sH_{21}(0)$. It is easy to verify that this is precisely the same approximation we would have obtained by using a forward Euler step to approximate the solution at $s$ of the differential equation satisfied by $Y$.*

2. *By our construction, we see that $Y_0$ in (33) is a second order approximation (in $s$) to the exact solution $Y(s)$. More precisely, from Taylor expansion at $s = 0$ of the exact solution $Y(s)$, we immediately get that $Y(s) = s\dot{Y}(0) + O(s^2) = sH_{21}(0) + O(s^2)$. On the other hand, the initial guess $Y_0 = 0$ is only an $O(s)$ approximation to $Y(s)$.*

3. *In practice, in (32), we do not have a close expression for $\dot{A}(0)$, but we can replace it by the difference quotient $(1/s)(A(s) - A(0))$. With this choice, the value of $Y_0$ (defined as in (33)) turns out to be the solution of the Sylvester equation*
   $$T_{22}(0)Y_0 - Y_0 T_{11}(0) = -E_{21}. \tag{34}$$
   *Since $\dot{A}(0) = (1/s)(A(s) - A(0)) + O(s^2)$, the value of $Y_0$ in (34) is still a $O(s^2)$ approximation to the exact $Y(s)$.*

[2] Use of the Euler step is an accepted standard in continuation algorithms, and it is the strategy implemented in `AUTO`.

4. *The estimates (and proof) of convergence for (26) and (27) relied on the value of $\kappa$ in (28) to be sufficiently small. Much the same arguments can be used for the refined guess in (33), by the following simple modification. In $U(s)$ in (21), we let $Y(s) = Y_0 + \Delta Y$ and use this form of $U$ in the invariant subspace relation (23). So doing, we obtain the new Riccati equation*

$$(\widehat{T}_{22} - Y_0\widehat{T}_{12})\Delta Y - \Delta Y(\widehat{T}_{11} + \widehat{T}_{12}Y_0) = -F(Y_0) + \Delta Y T_{12}\Delta Y. \tag{35}$$

*Now, it is a simple verification that iterative solution of (35) by either of the two iterations (26) or (27), appropriately reformulated for the unknown $\Delta Y$, produces precisely the same sequence as the original iterations on the unknown $Y$ had they been started with initial guess $Y_0$ in (33). As a consequence, the convergence results we quoted after (28) hold unchanged, except that we have a new $\kappa$ value, call it $\tilde{\kappa}$ :*

$$\tilde{\kappa} = \frac{\|\widehat{T}_{12}\|_F \|F(Y_0)\|_F}{sep^2(\tilde{T}_{11}, \tilde{T}_{22})}\ , \quad \tilde{T}_{11} = \widehat{T}_{11} + \widehat{T}_{12}Y_0, \;\; \tilde{T}_{22} = \widehat{T}_{22} - Y_0\widehat{T}_{12}\,. \tag{36}$$

*We now look at the ratio*

$$\frac{\tilde{\kappa}}{\kappa} \;=\; \frac{\|F(Y_0)\|_F}{\|E_{21}\|_F}\left(\frac{sep(\widehat{T}_{11}, \widehat{T}_{22})}{sep(\tilde{T}_{11}, \tilde{T}_{22})}\right)^2 .$$

*By our previous remark, $F(Y_0)$ is close to $0$ at second order in $s$, whereas $E_{21}$ is only first order close to $0$. To compare $sep(\tilde{T}_{11}, \tilde{T}_{22})$ with $sep(\widehat{T}_{11}, \widehat{T}_{22})$, we realize that $sep(\tilde{T}_{11}, \tilde{T}_{22})$ is a second order approximation to $sep(T_{11}(s), T_{22}(s))$, whereas $sep(\widehat{T}_{11}, \widehat{T}_{22})$ is only a first order approximation to the same quantity. Therefore, in first approximation, $\left(\frac{sep(\widehat{T}_{11}, \widehat{T}_{22})}{sep(\tilde{T}_{11}, \tilde{T}_{22})}\right)^2$ is $1 + O(s/sep(T_{11}(s), T_{22}(s)))$. To summarize, we have that*

$$\tilde{\kappa}/\kappa = O(s), \tag{37}$$

*with an obvious improvement in the radius of convergence for both iterations (26) and (27).*

5. *By creating the initial guess $Y_0$ from (33) we need to solve the Sylvester equation (32), or (34). Typically, this involves Schur reduction of the matrices $T_{11}(0)$ and $T_{22}(0)$, which is somewhat expensive. However, suppose we use the Newton iteration (27); then, at each iteration step we need to Schur factor the matrices $\widehat{T}_{11} + Y_{k-1}\widehat{T}_{12}$ and $\widehat{T}_{22} - \widehat{T}_{12}Y$. At convergence of Newton's method, these are precisely the matrices we will need for the next continuation step in (32). Therefore, no extra factorizations are needed in this case. For a practical comparison of the iterations for the two initial guesses (33) and $Y_0 = 0$, we refer to the next section.*

**Remark 7** *An alternative to our approach for continuing orthogonal invariant subspaces could be easily obtained as follows. Let $Q(0)$ : $Q^T(0)A(0)Q(0)$ be a Schur factorization of $A(0)$, $Q(0) = [Q_1(0),\ Q_2(0)]$ as usual, so that (say) $Q_1(0)$ spans the invariant subspace relative to the eigenvalues with positive real parts. Next, consider $A(s)$ which has the same* *number as $A(0)$ of eigenvalues with positive real part. Then, one may think to take the ordered Schur factorization of the matrix $A(s)$ : $P^TA(s)P$, and partition $P = [P_1,\ P_2]$ so that $P_1$ spans the invariant subspace relative to the eigenvalues with positive real part of $A(s)$. In general, $P_1$ is not smooth (i.e., it is not true that $P_1 = Q_1(0) + s\dot{Q}_1(0) + \dots$). To enforce smoothness, we may solve an orthogonal Procrustes problem and replace $P_1$ with $P_1V$, where $V$ is the orthogonal polar factor of $P_1^TQ_1(0)$ (see [12]). Essentially, this is the approach used by Beyn, [2], and then implemented in* `HomCont`*, [4]. However, we believe that the approach we have adopted is preferable to the one we just outlined. For one thing, the latter is usually more expensive than our approach. Moreover, unlike the one we just outlined, our approach to continuation of invariant subspaces is more in tune with the original continuation problem: continuation of the subspaces influences the continuation step, and using first derivative information (as we did to obtain (33)) is bound to reflect genuine difficulties of the original differential equation (such as nearing a bifurcation) into the continuation algorithm.*

# 4 Example: Heteroclinic orbits in a 4-D singular perturbation problem.

Consider the problem of finding traveling wave front solutions to the FitzHugh-Nagumo equations with two diffusive variables

$$v_t = v_{xx} + v(v-a)(1-v) - w, \quad w_t = \delta w_{xx} + \epsilon(v - \gamma w), \tag{38}$$

for small positive $\epsilon$, while $\delta$ ranging between a small and large value. For $\delta$ small this is a singularly perturbed reaction-diffusion system. In moving coordinates, $v_1 = v(z),\ v_2 = v'(z)$, $w_1 = w(z)$, $w_2 = w'(z)$ with $z = t + cx$, the reduced ODE is

$$\begin{aligned} v_1' &= v_2, \\ v_2' &= cv_2 - v_1(1-v_1)(v_1-a) + w, \\ w_1' &= w_2, \\ w_2' &= [cw_2 - \epsilon(v_1 - \gamma w_1)]/\delta. \end{aligned} \tag{39}$$

In [9] we located a heteroclinic orbit of (39). Here we first reproduced this result with our new code and then located a heteroclinic orbit with $\delta = \epsilon = 0.001$, $\gamma = 13.23529$, $a = 0.3$, and $c = 0.2571271$. In this case $n_0 = n_1 = 2$, where the relevant eigenvalues are $\mu_{0,1} = 0.6958$, $\mu_{0,2} = 257.2$, $\mu_{1,1} = -0.4247$, and $\mu_{1,2} = -0.06553$. We then performed a two parameter continuation in $(\delta, c)$ in the direction of increasing $\delta$, see Fig. 1 for the bifurcation diagram and Fig. 2 for some typical solutions in $(v_1, t)$ coordinates.

We have implemented 4 methods to solve the Riccati equation: (i) the "Simple Iteration" (26) with zero initial guess, (ii) Newton's method (27) with zero initial guess, (iii) Simple Iteration (26) with Euler initial guess, and (iv) Newton's method with Euler initial guess. The numerical results agree with our theoretical results. Namely, the choice of the Euler initial guess is a big improvement with respect to the simpler zero guess. Specifically, for the same pseudo-arclength continuation step (it took 433 such steps to compute the branch), (iii) and (iv) required much fewer iterations to converge than (i) and (ii). Some insight in the advantages gained by use of the Euler guess is obtained by considering the typical convergence behavior of the four methods (i)–(iv): on average, method (i) required 7 iterations for convergence, method (ii) needed 5, method (iii) needed slightly more than 4, and method (iv) required less than 3 iterations for convergence. In one exceptional case, methods (i) and (ii) required as many as 25 iterations. Finally, (i) and (ii) failed to converge in some cases (at the end of the continuation), because the continuation step was not small enough, whereas (iii) and (iv) never failed to converge.

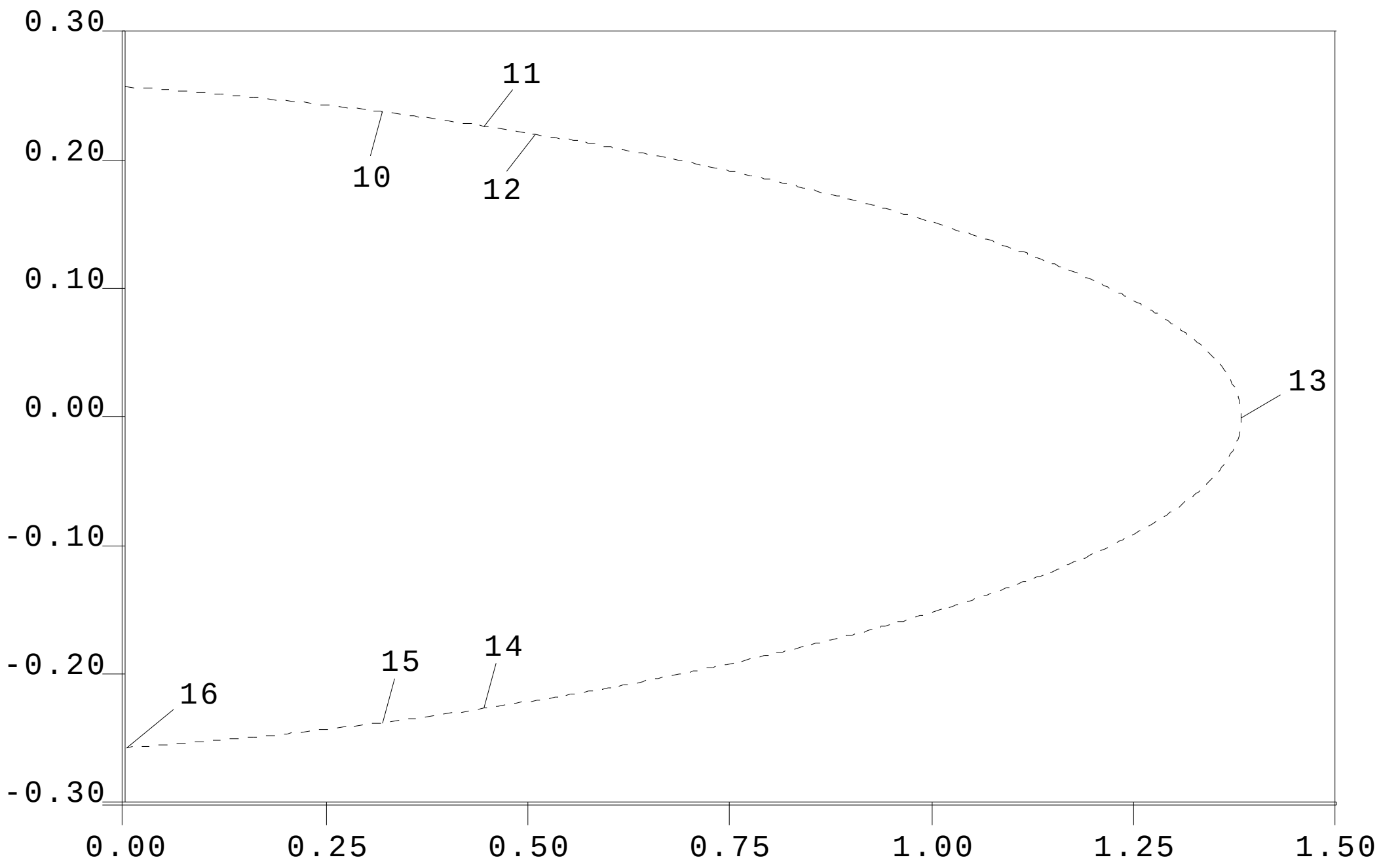

Figure 1. *Bifurcation diagram in* $(\delta, c)$ *coordinates. At label* 10 $(\delta, c) = (0.3198, 0.2376)$*, the eigenvalues are:* $\mu_{0,1} = \mu_{0,2} = 0.7406$*,* $\mu_{1,1} = -0.4357$*, and* $\mu_{1,2} = -0.06502$*. At label* 11 $(\delta, c) = (0.4462, 0.2265)$*, the eigenvalues are:* $\mu_{0,1} = \mu_{0,2} = 0.6204$*,* $\mu_{1,1} = -0.4409$*, and* $\mu_{1,2} = -0.06565$*. The eigenvalues are complex between the labels* 10 *and* 11*. At label* 12 $(\delta, c) = (0.5085, 0.2202)$*, the eigenvalues are:* $\mu_{0,1} = 0.5098$*,* $\mu_{0,2} = 0.6538$*,* $\mu_{1,1} = -0.4436$*, and* $\mu_{1,2} = -0.06612$*, where* $|\mu_{1,1}| + |\mu_{1,2}| = \mu_{0,1}$*. At label* 13 $(\delta, c) = (1.383, 0)$*, the eigenvalues are:* $\mu_{0,1} = -\mu_{1,2} = 0.1099$*,* $\mu_{0,2} = -\mu_{1,1} = 0.5454$*. And the part of the branch below the* $\delta$*-axis is symmetric to the one above the* $\delta$*-axis.*

Labels 10, 11, 12, 14, and 15 mark local bifurcations (of the eigenvalues), while label 13 marks a global bifurcation (the detailed analysis of these bifurcations will be given elsewhere). For comparison, we repeated the computation of the above branch using `HomCont` (in `AUTO97`), see Fig. 3. Note that `HomCont` could not continue the original branch beyond the global bifurcation point and instead switched the branches. We also note that by varying the continuation step size in our code we could select the desired branch at Label 13, while it was not possible to achieve this with `HomCont`. This confirms our theoretical insight that our numerical method is more in tune with the original continuation problem.

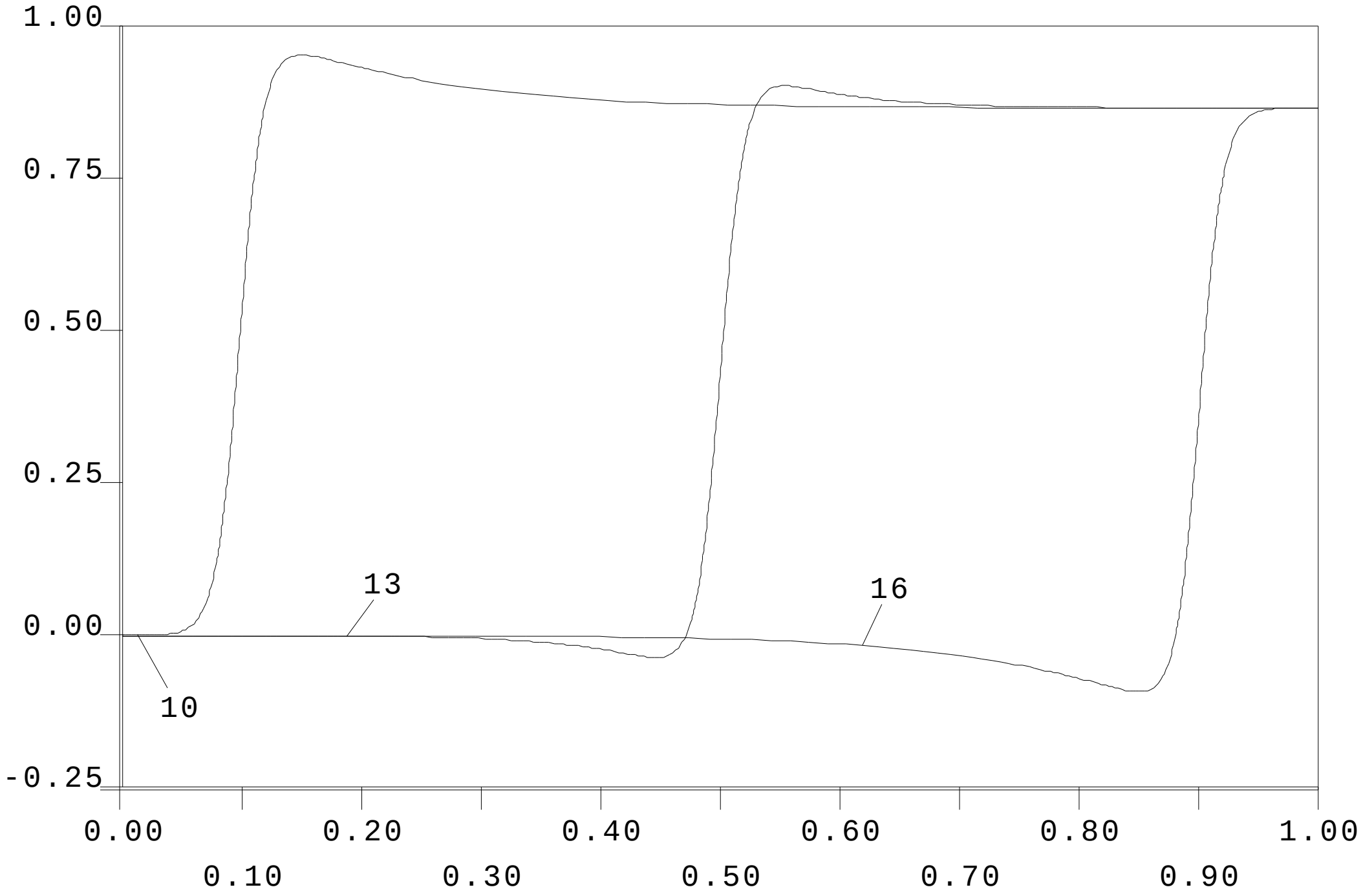

Figure 2. *Some typical solutions in* $(v_1, t)$ *coordinates.*

**Acknowledgments.** The authors wish to thank Z. Bai, University of Kentucky, for helpful discussions and for providing a subroutine for solving a Sylvester equation and W. -J. Beyn, University of Bielefeld, for helpful comments on an earlier version of the manuscript. The third author is also grateful to the IMA of University of Minnesota, to the Department of Mathematics of the University of Utah, to the Computer Science Division of the University of California at Berkeley, and to the School of Mathematics of Georgia Tech, for their hospitality during parts of his sabbatical stay, when this work was completed.

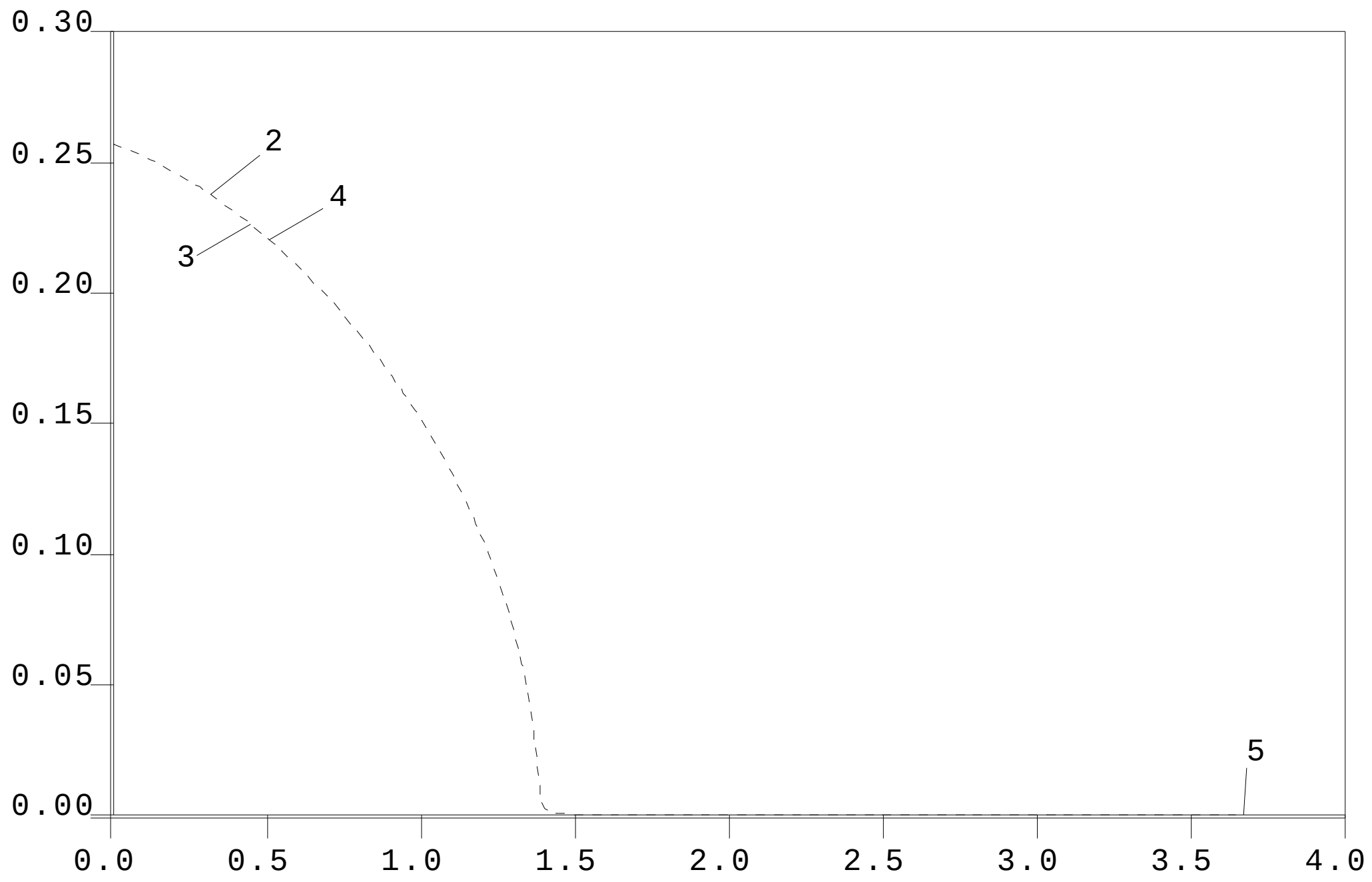

Figure 3. *Bifurcation diagram in $(\delta, c)$ coordinates computed by* `HomCont` *(AUTO97).*